\documentclass[a4paper]{amsart}
\usepackage{amsmath,amssymb,amsthm}
\usepackage[margin=1in]{geometry}
\usepackage{xcolor}
\usepackage{enumerate}
\usepackage{comment}
\usepackage[dvipsnames]{xcolor}
\usepackage[colorlinks, citecolor = blue]{hyperref}
\newtheorem{theorem}{Theorem}[section]
\newtheorem{lemma}[theorem]{Lemma}
\newtheorem{proposition}[theorem]{Proposition}

\theoremstyle{definition}
\newtheorem{definition}[theorem]{Definition}
\newtheorem{remark}[theorem]{Remark}
\usepackage{todonotes}

\numberwithin{equation}{section}

\newcommand{\R}{\mathbb{R}}

\newcommand{\1}{\mathbf{1}}
\newcommand{\avg}[2]{\langle #1\rangle_{#2}}
\newcommand{\D}{\mathcal{D}}
\DeclareMathOperator*{\essinf}{ess\,inf}
\DeclareMathOperator*{\esssup}{ess\,sup}
\newcommand{\pp}{p(\cdot)}
\newcommand{\ppr}{p'(\cdot)}
\newcommand{\qq}{q(\cdot)}
\newcommand{\A}{\mathcal{A}}
\newcommand{\Sp}{\mathcal{S}}

\author{Francisco Gon\c{c}alves}
\address{Delft Institute of Applied Mathematics, Delft University of
Technology, P.O. Box 5031, 2600 GA Delft, The Netherlands}
\email{f.g.j.diasdecarvalho@tudelft.nl}
\author{Tuomas Oikari}
\address{Aalto University, Department of Mathematics and Systems Analysis, P.O. Box 11100,
	FI-00076 Aalto, Finland}
\email{tuomas.oikari@aalto.fi}
\thanks{
T.O. was supported by the Research Council of Finland, projects 
360184 and 371637. He is a member of the Finnish Centre of Excellence in Randomness and Structures.}

\author{Troy Roberts}
\address{Washington University in St. Louis, Department of Mathematics, One Brookings Drive, St. Louis, MO, 63130, USA}
\email{r.troy@wustl.edu}
\thanks{F.G. was supported by the Dutch Research Council (NWO) through grant \href{https://www.nwo.nl/en/projects/vividi223019}{VI.Vidi.223.019}.}
\title[Boundedness of the maximal function on weighted variable Lebesgue spaces]{An $A_1$ characterization of the boundedness of the maximal function on weighted variable Lebesgue spaces}

\date{}

\begin{document}

\begin{abstract}
We characterize the boundedness of the Hardy–Littlewood maximal operator $M$ on weighted variable Lebesgue spaces for exponents satisfying $1<p_-\le p_+<\infty$.
We show that $M$ is bounded on $L^{p(\cdot)}_w$ if and only if the pair $(p(\cdot),w)$ satisfies a variable Muckenhoupt condition and there exists $h\in L^1$ such that $w^{-1}h^{1/p(\cdot)}$ and $wh^{1/p'(\cdot)}$ belong to $A_1$. This characterization is new even in the unweighted setting.
To prove sufficiency, we use sparse domination together with a local modular averaging estimate whose error terms are summable over sparse families to \(\|h\|_{L^1}\). Conversely, assuming $M$ is bounded, we construct $h$ using a simultaneous Rubio de Francia algorithm. To motivate this new characterization, we show that the strong Muckenhoupt condition, which characterizes the boundedness of $M$ in the unweighted case by a theorem of Diening, is not sufficient in the weighted setting. As applications, we give short proofs and extensions of known sufficient conditions.
\end{abstract}
\maketitle

\section{Introduction and main results}
The Hardy-Littlewood maximal operator \cite{HL30} is defined for $f\in L^1_{\mathrm{loc}}(\mathbb{R}^n)$ by
\[
Mf(x):=\sup_{Q\ni x}\langle |f|\rangle_Q,
\]
where $\langle f\rangle_Q:=\frac{1}{|Q|}\int_Q f(x)\mathrm d x$ and the supremum is taken over all cubes containing the point $x$.
For $p\in(1,\infty)$ and a weight $w>0$, the weighted Lebesgue space $L^p_w = L^p_w(\mathbb{R}^n)$ consists of those functions for which the norm $\|f\|_{L^p_w(\mathbb{R}^n)}:=\|fw\|_{L^p(\mathbb R^n)}$ is finite. We now mainly drop the ambient space $\mathbb{R}^n$ from notation.
Muckenhoupt \cite{Muc72} proved that $M$ is bounded on $L^p_w$ if and only if the weight $w$ satisfies the $A_p$ condition,
\[
[w]_{A_p} := \sup_Q |Q|^{-1}
{\|\mathbf 1_Q\|_{L^p_w}
      \|\mathbf 1_Q\|_{L^{p'}_{w^{-1}}}}
<\infty,
\]
where $p'=p/(p-1)$.
Moreover, this $A_p$ Muckenhoupt condition characterizes the boundedness of all Calder\'on--Zygmund operators on $L^p_w$, as shown by Coifman and Fefferman \cite{CF74}.

A variable exponent Lebesgue space $L^{p(\cdot)}(\R^n)$  allows the integrability to change throughout the domain through a variable exponent $p(\cdot):\R^n\to(1,\infty)$.
 Throughout the article we assume that $p(\cdot)$ is measurable and that 
\begin{equation}\label{eq:bounds}
1<p_-:=\essinf_{x\in\R^n}p(x)
\le \esssup_{x\in\R^n}p(x)=:p_+<\infty.
\end{equation}
The space $L^{p(\cdot)}(\R^n)$ consists of the measurable functions with finite Luxemburg norm
\[
\|f\|_{L^{p(\cdot)}(\R^n)}
:=\inf\biggl\{\lambda>0:
\int_{\R^n}\biggl|\frac{f(x)}{\lambda}\biggr|^{p(x)}
\,\mathrm d x\le1\biggr\}.
\]

Beyond their intrinsic interest, variable exponent Lebesgue spaces have been studied for their applications to PDEs and variational problems with non-standard growth \cite{AcMi02,DienRuz03,Fan03PxLap,HarHasLeNuo10,LLP10,Mingione06,Ruz00}.
One of the main problems in the field has been to characterize exactly when $M:L^{p(\cdot)}\to L^{p(\cdot)}$ maps boundedly. 

The variable-exponent analogue of the classical Muckenhoupt condition is the $A_{\pp}$ condition. We say that $1 \in A_{\pp}$ if 
\[
[1]_{A_{\pp}} := \sup_Q
|Q|^{-1} \|\mathbf 1_Q\|_{L^{p(\cdot)}}
      \|\mathbf 1_Q\|_{L^{p'(\cdot)}}
<\infty,
\]
where we write $p'(\cdot)=p(\cdot)/(p(\cdot)-1)$ for the conjugate exponent. 
It is equivalent to uniform boundedness of the averaging operators $f\mapsto\langle f\rangle_Q1_Q$ on $L^{\pp}$, and hence is necessary for boundedness of $M$ on $L^{p(\cdot)}$.
Kopaliani \cite{Kop07} proved that $1 \in A_{\pp}$ is sufficient when $p(\cdot)$ is constant outside a ball.
However, his subsequent counterexample \cite{Kop08} shows that the $A_{\pp}$ condition alone does not suffice for general exponents, which is in sharp contrast to Muckenhoupt's characterization for $L^p_w$.

Diening \cite{Die04} proved that a uniform local log-H\"older condition on $p(\cdot)$ suffices for boundedness of $M$ on $L^{p(\cdot)}$, provided that $p(\cdot)$ is constant outside a ball.
Cruz-Uribe, Fiorenza and Neugebauer \cite{CFN03,CFN04} relaxed this assumption at infinity to log-H\"older decay to a constant, yielding the familiar sufficient condition of global log-H\"older continuity.
Keeping the same local hypothesis, Nekvinda \cite{Nek04} weakened the requirement at infinity further by replacing it with an integral condition.
Continuity is not necessary, however, as Lerner \cite{Ler05} showed by constructing discontinuous exponents for which $M$ is bounded.

The first necessary and sufficient condition was given by Diening, who characterized boundedness of $M$ by the uniform boundedness of sums of averaging operators over pairwise disjoint cubes \cite[Theorem~8.1]{Die05}.
Using Diening's criterion, Lerner \cite[Theorem~1.12]{LerCrit26} obtained a characterization directly in terms of $p(\cdot)$, combining the $A_{\pp}$ condition with an integral condition at infinity.
More recently, Lerner \cite[Theorem~1.3]{LerAinf26} characterized boundedness by an $A_\infty$-type condition for both $p(\cdot)$ and $p'(\cdot)$, comparing norms of step functions on disjoint cubes before and after restriction to large subsets.
Both criteria rely on Diening's criterion and require uniform estimates over arbitrary families of pairwise disjoint cubes, which can make them difficult to verify for a given exponent.

Our first result provides a more easily verifiable characterization of the boundedness of $M$.
In addition to assuming $1 \in A_{\pp}$, we require the existence of a positive integrable function $h$ such that two variable powers of $h$ belong to $A_1$.
We recall that a locally integrable function $u>0$ belongs to $A_1$ (see \cite{Muc72}) if
\begin{align}\label{eq:defnA1}
    [u]_{A_1}:=\sup_Q
\frac{\avg{u}{Q}}{\essinf_{x\in Q}u(x)}<\infty.
\end{align}
Our first result is as follows. 
\begin{theorem}\label{thm:main-unweighted}
Let $p(\cdot)$ satisfy \eqref{eq:bounds}.
Then $M$ is bounded on $L^{p(\cdot)}$ if and only if $1 \in A_{\pp}$ and there exists a positive $h\in L^1$ such that
\begin{equation}\label{eq:compatible-unweighted}
h^{1/\pp}\in A_1
\qquad\text{and}\qquad
h^{1/\ppr}\in A_1.
\end{equation}
\end{theorem}

We call such a function $h$ a \emph{compatible density} for $p(\cdot)$.
As in the classical weighted setting, these conditions also characterize boundedness of all Calder\'on--Zygmund operators on $L^{p(\cdot)}$.
Indeed, Diening \cite[Theorems~8.1 and~8.14]{Die05} proved that boundedness of $M$ implies boundedness of these operators, while the converse follows from Rutsky's theorem \cite[Theorem~2]{Rut16}.

We also study the problem of when $M$ is bounded on the weighted variable exponent Lebesgue spaces.
For a weight $w>0$, we recall that $L^{p(\cdot)}_w$ consists of measurable functions $f$ such that the norm
\[
\|f\|_{L^{p(\cdot)}_w}
:=\|fw\|_{L^{p(\cdot)}}
\]
is finite.
We say $w \in A_{\pp}$ if 
\[
[w]_{A_{\pp}}:=\sup_Q
|Q|^{-1} \|\mathbf 1_Q\|_{L^{p(\cdot)}_w}
      \|\mathbf 1_Q\|_{L^{p'(\cdot)}_{w^{-1}}} < \infty.
\]
We remark that for $p(\cdot) = p$ constant this recovers the classical $A_p$ condition of Muckenhoupt. Moreover, if $\pp$ satisfies \eqref{eq:bounds} and $w \in A_{\pp}$, then $w^{\pp}, w^{-\ppr} \in L^1_{\mathrm{loc}}$.

Kokilashvili and Samko \cite[Theorem~A]{KS04} characterized boundedness of $M$ on $L^{p(\cdot)}_w$ for power weights $w(x)=|x-x_0|^\beta$ on bounded domains containing $x_0$, assuming local log-H\"older continuity of $p(\cdot)$.
For general weights on $\R^n$, Cruz-Uribe, Diening and H\"ast\"o \cite[Theorem~1.3]{CDH11} proved that $w \in A_{\pp}$ characterizes boundedness of $M$ on $L^{p(\cdot)}_w$ when $p(\cdot)$ is globally log-H\"older continuous.
Cruz-Uribe, Fiorenza and Neugebauer \cite[Theorem~1.5]{CFN12} obtained the same characterization by a different argument. By taking $w = 1$, one sees that these continuity assumptions are not necessary by Lerner's example in the unweighted setting \cite{Ler05}. Furthermore, as we show in Proposition \ref{prop:dc-counterexample}, Diening's characterization of $M : L^{\pp} \to L^{\pp}$ by uniform boundedness of sums of averaging operators over pairwise disjoint cubes cannot be extended to the weighted setting. This prompts the natural question: how can we characterize the pairs $(\pp,w)$ for which $M : L^{\pp}_w \to L^{\pp}_w$? 

The following Theorem \ref{thm:main} is the main contribution of the paper. It recognizes that a weighted version of the compatible density condition together with the $A_{\pp}$ condition provides a full characterization (provided \eqref{eq:bounds} holds) of the boundedness of the maximal function on weighted variable exponent Lebesgue spaces. In particular, this result does not assume log-H\"older continuity of the exponent, and hence is the first characterization of boundedness for a rough exponent $p(\cdot)$ in the weighted setting without any a priori continuity assumptions.
\begin{theorem}\label{thm:main}
Let $p(\cdot)$ satisfy \eqref{eq:bounds} and let $w$ be a weight.
Then $M$ is bounded on $L^{p(\cdot)}_w$ 
if and only if
$w \in A_{\pp}$ and there exists a positive function $h\in L^1$ such that
\begin{equation}\label{eq:compatible}
w^{-1}h^{1/p(\cdot)}\in A_1
\qquad\text{and}\qquad
wh^{1/p'(\cdot)}\in A_1.
\end{equation}
\end{theorem}
We call such a function $h$ a \emph{compatible density} for $(p(\cdot),w)$. One way to interpret the existence of a compatible density is through the lens of Jones factorization. It is a classical result that, in the multiplier notation, a weight $w \in A_p$ if and only if there exist $u,v \in A_1$ such that $w = u^{1/p}v^{-1/p'}$ \cite{Jon80}. Taking $u = wh^{1/\ppr}$ and $v = w^{-1}h^{1/\pp}$, we see that a compatible density gives us the decomposition $w = u^{1/\pp}v^{-1/\ppr}$ for $u,v \in A_1$ and $uv \in L^1$.
For a constant exponent $p$, the fact that the $A_1$ weights $u,v$ can be chosen so that $uv \in L^1$ is implicit in the proof given in \cite{CJR83}.  

We note that the two $A_1$ weights in \eqref{eq:compatible} are interchanged under mapping $(p(\cdot),w)\mapsto (p'(\cdot),w^{-1})$ (and back) and also
$[w]_{A_{\pp}}=[w^{-1}]_{A_{\ppr}}$. Thus, the same criterion characterizes boundedness of $M$ on $L^{p'(\cdot)}_{w^{-1}}$.
The equivalence of these two maximal bounds is due to Lerner \cite[Corollary~1.5]{Ler26}, which extended his earlier result where the assumption $w^{p(\cdot)}\in A_\infty$ was required \cite[Theorem~1.1]{Ler17}.
We use \cite[Corollary~1.5]{Ler26} in the proof of necessity when constructing a compatible density $h$ for $(p(\cdot),w).$

\subsection*{Outline of the paper} 
Section \ref{sec:Preliminaries} is for preliminaries.

Section \ref{sec:Core} is the core of the paper and contains the proof of Theorem \ref{thm:main}. 

Section \ref{sec:proof1} contains the sufficiency of $w \in A_{\pp}$ and the existence of a compatible density for the boundedness of the maximal function.
Pointwise sparse domination (Proposition \ref{prop:prelim-sparse-domination}) reduces the problem to proving uniform bounds for sparse operators (Theorem \ref{thm:sparse}).
Given a sparse collection of cubes $\mathcal{S}$, let $E_Q$ denote their pairwise disjoint major subsets. 
The key technical result (Lemma \ref{lem:KEYlemma}) uses the $w \in A_{\pp}$ condition and the compatible density $h$ to control the weighted modular of the average $\langle f1_{E_Q}\rangle_Q$ by that of $f1_{E_Q}$ with an error proportional to $\|h1_{E_Q}\|_{L^1}.$ Since $E_Q$ are disjoint these errors are summable to $\|h\|_{L^1}$ and this will result in uniform bounds for sparse averaging operators.

Section \ref{sec:proof2} contains the necessity of $w \in A_{\pp}$ and the existence of a compatible density for the boundedness of the maximal operator. The condition $w \in A_{\pp}$ is immediate.
The construction of a compatible density is a nonlinear iteration based on the factorization arguments of Rubio de Francia \cite{RdF82} and Coifman, Jones and Rubio de Francia
\cite[proof of Theorem~II]{CJR83}.  
For the construction algorithm to work, we use Lerner's duality theorem \cite[Corollary~1.5]{Ler26} to obtain the boundedness of $M$ on the associate space $L^{p'(\cdot)}_{w^{-1}}.$
The boundedness of $M$ on $L^{p(\cdot)}_{w}$ and on $L^{p'(\cdot)}_{w^{-1}}$ keeps the increasing recursive sequence uniformly bounded in $L^1$, and its limit is a positive integrable function satisfying both $A_1$ conditions in \eqref{eq:compatible}.

Section \ref{sec:Applications} contains applications and demonstrates how the concept of compatible densities can be used to give concise proofs of known results, such as the sufficiency of $1 \in A_{\pp}$ and log-H\"older continuity at infinity for $M : L^{\pp} \to L^{\pp}$.

Finally, in Section \ref{sec:failure-diening}, we give an example showing that Diening's criterion of uniform boundedness of sums of averaging operators over pairwise disjoint cubes cannot be extended to the weighted variable Lebesgue spaces.

\section{Preliminaries}\label{sec:Preliminaries}
In this section we recall basic notation and norm and modular properties of weighted variable Lebesgue spaces used throughout the paper. 

\subsection*{Basic definitions and notation}
We recall basic definitions and notational conventions.
    \begin{itemize}
    \item In this article a variable exponent is a measurable function $p(\cdot):\mathbb{R}^n\to (1,\infty)$. 
    \item We mostly suppress the ambient space $\mathbb{R}^n$.
  \item 
  Given
 a variable exponent $p(\cdot)$ and a weight
  $w>0$, the weighted variable exponent Lebesgue space $L^{p(\cdot)}_w$ consists
  of those measurable functions $f$ for which the norm
\begin{equation*}
  \|f\|_{p(\cdot),w} := \|f\|_{L^{p(\cdot)}_w}:=\inf\{\lambda>0:\rho_{p(\cdot),w}(f/\lambda)\leq 1\}
\end{equation*}
is finite, and
where the weighted modular is 
\begin{equation*}
	\rho_{p(\cdot),w}(f):=
	\int_{\mathbb R^n} |f(x)w(x)|^{p(x)}\,\mathrm d x.
		\end{equation*}
\item The conjugate exponent $p'(\cdot)$ is given by the relation $\frac{1}{p(\cdot)}+ \frac{1}{p'(\cdot)}=1.$
\end{itemize}
We use the following notation for absorbing constants.
\begin{itemize}
\item We write $A\lesssim B$ if there exists a constant $C>0$ such that $A\leq CB$. We write $A\simeq B$ if $A\lesssim B$ and $B\lesssim A.$ 
Subscripted or parenthesized data on constants $(C_{a,b}, C(a,b))$ and quantifiers $(\lesssim_{a,b}, \sim_{a,b})$ signify their dependence on that data. 
\end{itemize}

\subsection*{Basic properties of variable exponent Lebesgue spaces}
As the multiplication map $f\mapsto wf$ is an isometry from $L^{\pp}_w$ onto
$L^{\pp}$, the cited unweighted results automatically give the recorded weighted norm
and modular properties below. We first record the lattice and
monotone convergence properties. 
\begin{proposition}[{\cite[Theorem~3.2.13]{DHHR11}}]
\label{prop:Lp=Fatou}
The space $L^{\pp}_w$ is a Banach space. If $g\in L^{\pp}_w$ and
$f$ is measurable with $|f|\le|g|$ almost everywhere, then
$f\in L^{\pp}_w$ and $\|f\|_{\pp,w}\le\|g\|_{\pp,w}$.

Moreover, the norm has the Fatou property: if
$0\le f_k\uparrow f$ almost everywhere, then
\begin{equation}\label{eq:prelim-fatou}
\|f_k\|_{\pp,w}\uparrow\|f\|_{\pp,w}.
\end{equation}
\end{proposition}
The next norm--modular relations allow us to normalize functions to have modular one.

\begin{lemma}[{\cite[Lemma~3.4.2]{DHHR11}}]
\label{lem:prelim-norm-modular}
For all $f$, we have 
\begin{equation}\label{eq:prelim-unit-ball}
\rho_{p(\cdot),w}(f)\le1
\quad\Longleftrightarrow\quad
\|f\|_{\pp,w}\le1,
\end{equation}
and 
$$
\rho_{p(\cdot),w}(f)=1 \quad\Longleftrightarrow\quad\|f\|_{\pp,w}=1.
$$
In particular, if $0\ne f\in L^{p(\cdot)}_w$, then 
\begin{equation}\label{eq:prelim-normalized-modular}
\rho_{p(\cdot),w}\Bigl(\frac{f}{\|f\|_{\pp,w}}\Bigr)=1.
\end{equation}
\end{lemma}

In particular, for all cubes $Q$, by \eqref{eq:prelim-normalized-modular} such that $\1_Q\in L^{\pp}_w$, there holds that 
\begin{equation}\label{eq:prelim-indicator-normalization}
\int_Q
\biggl(\frac{w(x)}{\|\mathbf 1_Q\|_{\pp,w}}\biggr)^{p(x)}
\,\mathrm d x=1.
\end{equation}
Convex combinations can be controlled by the following modular estimate.
\begin{lemma}[{\cite[Propositions~2.7(4) and~2.14]{CUF13}}]
\label{lem:prelim-modular-estimates}
For $f,g\in L^{\pp}_w$ there holds that 
\begin{equation}\label{eq:prelim-modular-convexity}
\rho_{p(\cdot),w}\bigl(\theta f+(1-\theta)g\bigr)
\le\theta\rho_{p(\cdot),w}(f)+(1-\theta)\rho_{p(\cdot),w}(g),
\qquad 0\le\theta\le1.
\end{equation}
\end{lemma}
For pairings between $f\in L^{\pp}_w$ and $g\in L^{\ppr}_{w^{-1}}$, the
identity $fg=(wf)(w^{-1}g)$ gives the following form of H\"older's
inequality.
\begin{lemma}[{\cite[Lemma~3.2.20]{DHHR11}}]
\label{lem:Holder}
For $f\in L^{\pp}_w$ and $g\in L^{\ppr}_{w^{-1}}$ there holds that 
\begin{equation}\label{eq:prelim-holder}
\int_{\R^n}|f(x)g(x)|\,\mathrm d x
\le2\|f\|_{\pp,w}\|g\|_{\ppr,w^{-1}}.
\end{equation}
\end{lemma}
\begin{remark}\label{remark:LocalIntegrability}
    The assumptions $w^{p(\cdot)},w^{-p'(\cdot)}\in L^1_{\mathrm{loc}}$, which follow from $w \in A_{\pp}$ and \eqref{eq:bounds}, automatically give that $\|\mathbf 1_Q\|_{p(\cdot),w}<\infty$ and $\|\mathbf 1_Q\|_{p'(\cdot),w^{-1}}<\infty$, for all cubes $Q$. Applying Lemma~\ref{lem:Holder} thus
shows that $L^{p(\cdot)}_w\cup L^{p'(\cdot)}_{w^{-1}}\subset L^1_{\mathrm{loc}}$.
\end{remark}
The same pairing gives the following duality formula.

\begin{lemma}[{\cite[Corollary~3.2.14]{DHHR11}}]
\label{lem:dualityFormula}
For all $f\in L^{\pp}_w$ there holds that 
\begin{equation}\label{eq:prelim-duality}
\|f\|_{\pp,w}
\simeq
\sup_{\substack{\|g\|_{\ppr,w^{-1}}\le1}}
\int_{\R^n}f(x)g(x)\,\mathrm d x.
\end{equation}
\end{lemma}

\subsection{Preliminaries on the maximal function}
The main tool in proving sufficiency in Section \ref{sec:proof1} will be pointwise sparse domination of the maximal function by dyadic sparse operators. 
\begin{definition}\label{defn:DyadicSparse}  
We recall the following concepts.
\begin{itemize} 
   \item  A countable family of cubes $\mathcal{D}$ is \emph{dyadic}, if for all $Q,P\in\mathcal{D}$ there holds that $Q\cap P\in \{Q,P,\emptyset\}$.
   \item A countable family $\Sp$ of cubes is called \emph{sparse}, if there exist pairwise disjoint measurable sets $E_Q\subseteq Q$, $Q\in\Sp$ such that $|E_Q|\ge|Q|/2$. 
   \item If $\mathcal{S}$ is both sparse and dyadic, we say that $\mathcal{S}$ is a sparse collection of dyadic cubes. 
   \end{itemize}
\end{definition}
For a fixed sparse collection $\mathcal{S}$, we define the sparse averaging operators
\begin{equation}\label{eq:prelim-sparse-operators}
\A_{\Sp}f:=\sum_{Q\in\Sp}\avg{|f|}{Q}1_{E_Q}.
\end{equation}
We use the following pointwise sparse domination of the maximal function. 

\begin{proposition}[{\cite[Lemmas~2.3 and~2.6]{Ler17}}]
\label{prop:prelim-sparse-domination}
There exist dyadic lattices $\D_1,\ldots,\D_{3^n}$ such that for all $f\in L^{\infty}_c$ there exist dyadic sparse families $\Sp_j\subseteq\D_j$ such that
\begin{equation}\label{eq:prelim-sparse-domination}
Mf\lesssim_n\sum_{j=1}^{3^n}\A_{\Sp_j}f.
\end{equation}
\end{proposition}
The \emph{dyadic lattices} $\mathcal{D}_j$, $j=1,\dots,3^n$, are as in \cite[Section~2.3]{Ler17}, but their particular form is not relevant for us; after the domination we simply use that each $\mathcal{S}_j$ is a dyadic sparse collection (as per Definition \ref{defn:DyadicSparse}).
Uniform bounds for sparse operators on $L^{\pp}_w$ therefore
give boundedness of $M$ on $L^{\pp}_w$, and the estimate extends from
bounded compactly supported functions to all of $L^{p(\cdot)}_w$ by
applying it to $\min\{|f|,j\}1_{[-j,j]^n}$ and using
Proposition~\ref{prop:Lp=Fatou}.

For the necessity direction, to construct a compatible density in Section \ref{sec:proof2}, we use the following
duality property of the maximal function.
\begin{theorem}[{\cite[Corollary~1.5]{Ler26}}]\label{thm:prelim-maximal-duality}
Let $p(\cdot)$ satisfy \eqref{eq:bounds} and $w^{p(\cdot)},w^{-p'(\cdot)}\in L^1_{\mathrm{loc}}$. Then
\begin{equation}\label{eq:prelim-maximal-duality}
M:L^{\pp}_w\to L^{\pp}_w
\quad\Longleftrightarrow\quad
M:L^{\ppr}_{w^{-1}}\to L^{\ppr}_{w^{-1}}.
\end{equation}
\end{theorem}
For $w=1$, this result was obtained by Diening
\cite[Theorem~8.1]{Die05} through his characterization by uniform boundedness of sums of disjointly supported averaging operators. Lerner \cite[Theorem~1.1]{Ler17}
obtained the weighted result first under the additional assumption
$w^{p(\cdot)}\in A_\infty$ and then in \cite{Ler26} only under $w^{p(\cdot)},w^{-p'(\cdot)}\in L^1_{\mathrm{loc}}$.

\section{Proof of Theorem \ref{thm:main}}\label{sec:Core}
\subsection{Sufficiency}\label{sec:proof1}
We prove that the $[w]_{A_{p(\cdot)}}$ condition and the existence of a compatible density $h$ for the pair $(p(\cdot),w)$ are sufficient for the boundedness of the maximal function.
By Proposition~\ref{prop:prelim-sparse-domination}, it is
enough to prove these to be sufficient for uniform bounds for sparse operators as follows. 

\begin{theorem}\label{thm:sparse}
Let $p(\cdot)$ satisfy \eqref{eq:bounds} and let $w$ be a weight such that  $[w]_{A_{p(\cdot)}}<\infty$. Let $h$ be a compatible density for the pair $(p(\cdot), w).$
Let $\mathcal{S}$ be a countable sparse family of dyadic cubes. Then,
\[
\|\mathcal{A}_{\mathcal{S}}f\|_{p(\cdot),w}\lesssim \|f\|_{p(\cdot),w},
\]
for all $f\in L^{\infty}_c$, and
where the implicit constant depends on $p_-$, $[w]_{A_{p(\cdot)}}$, $\|h\|_{L^1}$,
  $[w^{-1}h^{\frac{1}{p(\cdot)}}]_{A_1}$, $
  [wh^{\frac{1}{p'(\cdot)}}]_{A_1}$.
\end{theorem}

The proof of Theorem~\ref{thm:sparse} uses the following local modular estimate, which can be viewed as a Jensen-type inequality in $L^{p'(\cdot)}_{w^{-1}}$ whose error term, controlled by $\int_E h$, is summable over sparse families. During the finalization process of our paper it came to our attention that a similar Jensen-type inequality, with error summable over disjoint families, was used by Shalukhina~\cite{Sha26} to give a short proof of Diening’s averaging criterion on spaces of homogeneous type.

\begin{lemma}\label{lem:KEYlemma} Let $p(\cdot)$ satisfy \eqref{eq:bounds} and let $w$ be a weight such that  $[w]_{A_{p(\cdot)}}<\infty$. Let $h$ be a compatible density for the pair $(p(\cdot), w).$ 
Let $Q$ be a cube and $E\subseteq Q$ such that $|E|\ge |Q|/2$.
Then, for all $f\ge0$ such that
$\int_E(f(x)/w(x))^{p'(x)}\,\mathrm d x\le1$ there holds that
\begin{equation}\label{eq:local}
\int_Q
\biggl(\frac{\avg{f1_E}{Q}}{w(x)}\biggr)^{p'(x)}\,\mathrm d x
\lesssim
\int_E\Bigl(\frac{f(x)}{w(x)}\Bigr)^{p'(x)}+h(x)
\,\mathrm d x,
\end{equation}
where the implicit constant depends on $p_-$, $[w]_{A_{p(\cdot)}}$, $\|h\|_{L^1}$, $
  [w^{-1}h^{1/p(\cdot)}]_{A_1}$, $
  [wh^{1/p'(\cdot)}]_{A_1}$.
\end{lemma}
\begin{remark}
    The variable exponent analogue of Jensen's inequality, $\int_Q \langle f \rangle_Q^{p(x)} \, dx \lesssim \int_Q f(x)^{p(x)} \, dx$, does not hold uniformly over all cubes. Jensen-type inequalities were used by Kopaliani \cite[Lemma 2.4]{Kop07} and Lerner \cite[Lemma 3.1]{Ler10}, where instead of a summable error, they impose restrictions on the size of the average and norm of $f$. The idea of a Jensen-type inequality with integrable error term was utilized in the proof that log-H\"older continuity is sufficient for $M : L^{\pp} \to L^{\pp}$ \cite[Corollary 4.2.5]{DHHR11}. Their summable error takes the form of an integrable power weight; we make this connection explicit through the proof of Theorem \ref{thm:log holder}.  
\end{remark}

Assuming Lemma~\ref{lem:KEYlemma}, we now prove
Theorem~\ref{thm:sparse}.
\begin{proof}[Proof of Theorem~\ref{thm:sparse}]
By Proposition \ref{prop:Lp=Fatou} we may assume that $\mathcal{S}$ is finite. 
After bounding $\mathcal{A}_{\mathcal{S}}f\leq \mathcal{A}_{\mathcal{S}}|f|$, we may assume that $f\geq 0.$ 
Then, by Lemma \ref{lem:dualityFormula},
\begin{align}\label{eq:proofSuff1}
     \|\mathcal{A}_{\mathcal{S}}f\|_{p(\cdot),w} \simeq
\sup_{\substack{g\ge0\\\|g\|_{p'(\cdot),w^{-1}}= 1}}
\int_{\R^n}\mathcal{A}_{\mathcal{S}}fg.
\end{align}
Fix such a $g$ as in the supremum on the right-hand side and let 
\[
\mathcal{A}_{\mathcal{S}}^*g := 
\sum_{Q\in\mathcal{S}}\avg{g1_{E_Q}}{Q}1_Q
\]
denote the adjoint of $\mathcal{A}_{\mathcal{S}}.$
By Tonelli's theorem and Lemma \ref{lem:Holder} we have 
\begin{align*}
  \int_{\R^n}\mathcal{A}_{\mathcal{S}}fg
=\sum_{Q\in\mathcal{S}}\avg{f}{Q}\int_{E_Q}g\,\mathrm d x
=\int_{\R^n}f
 \Bigl(
\sum_{Q\in\mathcal{S}}\avg{g1_{E_Q}}{Q}1_Q
\Bigr) &= \int f\mathcal{A}_{\mathcal{S}}^*g \lesssim \|f\|_{p(\cdot),w}\|\mathcal{A}_{\mathcal{S}}^*g\|_{p'(\cdot),w^{-1}}.
\end{align*}
Thus it remains to prove that 
\begin{equation}\label{eq:dualSparse1}
\|\mathcal{A}_{\mathcal{S}}^*g\|_{p'(\cdot),w^{-1}}\lesssim \|g\|_{p'(\cdot),w^{-1}}.
\end{equation}

Since $\mathcal{S}$ is finite, we have $\mathcal{A}_{\mathcal{S}}^*g\in L^{\infty}_{c}$. Thus, $\|\mathcal{A}_{\mathcal{S}}^*g\|_{p'(\cdot),w^{-1}} \leq \|\mathcal{A}_{\mathcal{S}}^*g\|_{\infty}\|1_{\mathrm{supp}(\mathcal{A}_{\mathcal{S}}^*g)}\|_{p'(\cdot),w^{-1}} < \infty,$ using the assumption $w^{-p'(\cdot)}\in L^1_{\mathrm{loc}}$, see Remark \ref{remark:LocalIntegrability}.
Hence, assuming $\mathcal{A}_{\mathcal{S}}^*g\not=0$ (if not, there is nothing to prove), we see that \eqref{eq:dualSparse1} is equivalent, after renormalization, to
\begin{align}\label{eq:dualSparse2}
    1= \|\mathcal{A}_{\mathcal{S}}^*\Tilde{g}\|_{p'(\cdot),w^{-1}}\lesssim \|\Tilde{g}\|_{p'(\cdot),w^{-1}},\quad\text{where}\quad \Tilde{g} := \frac{g}{\|\mathcal{A}_{\mathcal{S}}^*g\|_{p'(\cdot),w^{-1}}}.
\end{align}
To conclude the proof, we verify \eqref{eq:dualSparse2}.

Fix $x$ and order the (possibly empty) finite sequence (by $\mathcal{S}$ being finite) of cubes containing
$x$ as $Q_1\subsetneq\cdots\subsetneq Q_J$ and set
\[
S_j:=\sum_{\substack{P\in\mathcal{S}\\P\supseteq Q_j}}
\avg{\Tilde{g}1_{E_P}}{P},
\qquad j=1,\dots,J,
\qquad S_{J+1}:=0.
\]
Then $\mathcal{A}_{\mathcal{S}}^*\Tilde{g}(x)=S_1$ and
$S_j-S_{j+1}=\avg{\Tilde{g}1_{E_{Q_j}}}{Q_j}$.
Note the elementary bound $t^r-s^r\le r(t-s)t^{r-1}$, valid for $0\le s\le t$ and $r\geq 1$. (This bound follows immediately from the mean value theorem.) We will use this with the substitutions $t = S_j$ and $s = S_{j+1}$ and $r=p'(x)$. Thus,
\begin{multline}\label{eq:gaghog}
\bigl(\mathcal{A}_{\mathcal{S}}^*\Tilde{g}(x)\bigr)^{p'(x)} = S_1^{p'(x)} - S_{J+1}^{p'(x)}
=\sum_{j=1}^J\bigl(S_j^{p'(x)}-S_{j+1}^{p'(x)}\bigr)\\
\le p'(x)\sum_{j=1}^J
\avg{\Tilde{g}1_{E_{Q_j}}}{Q_j}S_j^{p'(x)-1} \lesssim_{p_-}\sum_{Q\in\mathcal{S}}\avg{\Tilde{g}1_{E_Q}}{Q}
\biggl(\sum_{\substack{P\in\mathcal{S}\\P\supseteq Q}}
\avg{\Tilde{g}1_{E_P}}{P}\biggr)^{p'(x)-1}1_Q(x),
\end{multline}
where we bounded
\begin{align}\label{eq:proofSuff3}
    p'(x) = \left(1-\frac{1}{p(x)}\right)^{-1} \leq \left(1-\frac{1}{p_-}\right)^{-1} \lesssim_{p_-} 1.
\end{align}
Multiplying both sides of \eqref{eq:gaghog} by $\|\Tilde{g}\|_{p'(\cdot),w^{-1}}^{-1}w^{-p'(\cdot)}$ and integrating gives 
\begin{multline}\label{eq:proofSuff4}
    \frac{1}{\|\Tilde{g}\|_{p'(\cdot),w^{-1}}} = \frac{1}{\|\Tilde{g}\|_{p'(\cdot),w^{-1}}}\int \bigl(\mathcal{A}_{\mathcal{S}}^*\Tilde{g}\bigr)^{p'(\cdot)}w^{-p'(\cdot)} \\ \lesssim_{p_-} \sum_{Q\in\mathcal{S}}\int_Q
\left\langle \frac{\Tilde{g}1_{E_Q}}{\|\Tilde{g}\|_{p'(\cdot),w^{-1}}}\right\rangle_Q
\biggl(\sum_{\substack{P\in\mathcal{S}\\P\supseteq Q}}
\avg{\Tilde{g}1_{E_P}}{P}\biggr)^{p'(\cdot)-1}
w^{-p'(\cdot)}.
\end{multline}
We next modify the integrand into a form compatible with Lemma \ref{lem:KEYlemma}.
Since $\sum_{\substack{P\in\mathcal{S}\\P\supseteq Q}}
\avg{\Tilde{g}1_{E_P}}{P}$ is constant on $Q$ and  $|E_Q|\geq \frac{1}{2}|Q|$ we have 
\begin{align*}
    \sum_{\substack{P\in\mathcal{S}\\P\supseteq Q}}
\avg{\Tilde{g}1_{E_P}}{P} = \left\langle  \sum_{\substack{P\in\mathcal{S}\\P\supseteq Q}}
\avg{\Tilde{g}1_{E_P}}{P}\right\rangle_{E_Q} \leq  2\left\langle  1_{E_Q}\sum_{\substack{P\in\mathcal{S}\\P\supseteq Q}}
\avg{\Tilde{g}1_{E_P}}{P}\right\rangle_{Q} \leq 2\left\langle 1_{E_Q}\mathcal{A}_{\mathcal{S}}^*\Tilde{g} \right\rangle_{Q}.
\end{align*}
Thus, denoting 
\[
G:=\frac{1}{2}\frac{\Tilde{g}}{\|\Tilde{g}\|_{p'(\cdot),w^{-1}}}+\frac{1}{2}\mathcal{A}_{\mathcal{S}}^*\Tilde{g}
\]
we have
\begin{align*}
    \mathrm{RHS}\eqref{eq:proofSuff4} &\leq \sum_{Q\in\mathcal{S}}\int_Q
\left\langle \frac{\Tilde{g}1_{E_Q}}{\|\Tilde{g}\|_{p'(\cdot),w^{-1}}}\right\rangle_Q
\biggl(2\left\langle 1_{E_Q}\mathcal{A}_{\mathcal{S}}^*\Tilde{g} \right\rangle_{Q}\biggr)^{p'(\cdot)-1}
w^{-p'(\cdot)} \\
&\leq \sum_{Q\in\mathcal{S}}\int_Q
\left\langle 2G1_{E_Q}\right\rangle_Q
\biggl(2\left\langle 2G1_{E_Q}\right\rangle_{Q}\biggr)^{p'(\cdot)-1}
w^{-p'(\cdot)} \lesssim_{p_-} \sum_{Q\in\mathcal{S}}\int_Q
\biggl(\frac{\avg{G1_{E_Q}}{Q}}{w}\biggr)^{p'(\cdot)},
\end{align*}
where in the last step we again used the bound \eqref{eq:proofSuff3}.
Consequently,
\begin{align}\label{eq:proofSuff5}
    \frac1{\|\Tilde{g}\|_{p'(\cdot),w^{-1}}} = \mathrm{LHS}\eqref{eq:proofSuff4} \lesssim_{p_-} \mathrm{RHS}\eqref{eq:proofSuff4}
\lesssim_{p_-}\sum_{Q\in\mathcal{S}}\int_Q
\biggl(\frac{\avg{G1_{E_Q}}{Q}}{w}\biggr)^{p'(\cdot)}.
\end{align}
To see that Lemma \ref{lem:KEYlemma} can be applied with $G$, 
by \eqref{eq:prelim-modular-convexity} we have 
\[
\rho_{p'(\cdot),w^{-1}}(G)
\le\frac12
\rho_{p'(\cdot),w^{-1}}\biggl(\frac{\Tilde{g}}{\|\Tilde{g}\|_{p'(\cdot),w^{-1}}}\biggr)
+\frac12\rho_{p'(\cdot),w^{-1}}(\mathcal{A}_{\mathcal{S}}^*\Tilde{g})=1.
\]
Thus, by
Lemma~\ref{lem:KEYlemma}, we continue the bound \eqref{eq:proofSuff5} as 
\begin{align*}
\frac1{\|\Tilde{g}\|_{p'(\cdot),w^{-1}}}
\lesssim_{p_-}\sum_{Q\in\mathcal{S}}\int_{E_Q}
\Bigl(\frac{G}{w}\Bigr)^{p'(\cdot)}+h
\le\rho_{p'(\cdot),w^{-1}}(G)+\|h\|_1\lesssim1.
\end{align*}
Thus $\|\Tilde{g}\|_{p'(\cdot),w^{-1}}\gtrsim1$, which proves
\eqref{eq:dualSparse2}.
\end{proof}

\subsection*{Proof of Lemma ~\ref{lem:KEYlemma}}

It remains to prove Lemma \ref{lem:KEYlemma}. We begin with an outline and the proof itself comes after that.

Fix a cube $Q$ and a modular budget $\rho_{p'(\cdot), w^{-1}}(f1_{E_Q}) =: t$. In a dream scenario we would like to
bound $\rho_{p'(\cdot), w}(\langle f1_{E_Q}\rangle_Q1_Q)\lesssim t$, for all $t\leq 1$ (and in this case the factor $\int_Eh$ would be absent from the right-hand side of \eqref{eq:local}). However, this bound only holds in a certain budget range $t\in [\epsilon_Q,1]$ (see \eqref{eq:defn:epsilonQ}); and below the budget (i.e. for $t<\epsilon_Q$) we use that $t$-budget falls into $\epsilon_Q$-budget. Thus,
\begin{align}\label{eq:heuristic}
      \rho_{p'(\cdot), w^{-1}}(\langle f1_{E_Q}\rangle_Q1_Q)\lesssim \begin{cases}
      \rho_{p'(\cdot), w^{-1}}(f1_{E_Q}),\quad  &\epsilon_Q\leq \rho_{p'(\cdot), w^{-1}}(f1_{E_Q}) \leq 1, \\ 
      \epsilon_Q,\quad  &0\leq \rho_{p'(\cdot), w^{-1}}(f1_{E_Q})< \epsilon_Q.
  \end{cases}
\end{align}
Since $\epsilon_Q\lesssim \int_Eh$ (directly from definition), \eqref{eq:heuristic} would yield \eqref{eq:local} and conclude the proof of Lemma \ref{lem:KEYlemma}.

The key tool to proving that a budget of $t\in[\epsilon_Q,1]$ gives the bound $\rho_{p'(\cdot), w^{-1}}(\langle f1_{E_Q}\rangle_Q1_Q)\lesssim t$ (Lemma \ref{lem: t adm modular inequality}) is that of $t$-admissible pairs (Definition \ref{defn:t-admissibility}), which is a sufficient condition for the bound to hold. 
Hence we would like to demonstrate the existence of $t$-admissible pairs for all budgets $t\in[\epsilon_Q,1].$  
When the budget is ``$t=1$", the $A_{p(\cdot)}$ condition provides a $1$-admissible pair (Lemma \ref{lem:pairs2}), while for ``$t=\epsilon_Q$" the compatible density provides an $\epsilon_Q$-admissible pair (Lemma \ref{lem:pairs2}). By interpolation
(Lemma \ref{lem: t adm geo mean}) we can construct from the endpoints $t=1$ and $t=\epsilon_Q$ admissible pairs for all budgets $t\in [\epsilon_Q, 1]$. This concludes the outline of the proof.

Then onto the proof itself.
For the rest of this section we denote 
\begin{align}
L:=\max\bigl\{1,
  [w^{-1}h^{1/p(\cdot)}]_{A_1},
  [wh^{1/p'(\cdot)}]_{A_1}\bigr\}.
\end{align}
\begin{definition}\label{defn:t-admissibility}
Let $1<r(\cdot)<\infty$ and $[w]_{A_{r(\cdot)}}<\infty$. For a cube $Q$ and constants
$a,b,t>0$, we call $(a,b)$ a \emph{$t$-admissible pair} on $Q$ for
$(r(\cdot),w)$ if
\begin{equation}\label{t-admissible}
\begin{aligned}
&\text{(i) }\int_Q\bigl(aw(x)\bigr)^{r(x)}\,\mathrm d x\le Lt,
\qquad
\text{(ii) }\int_Q\bigl(\frac{b}{w(x)}\bigr)^{r'(x)}\,\mathrm d x\le Lt,\\
&\text{(iii) }ab\ge\frac{t}{2[w]_{A_{r(\cdot)}}|Q|}.
\end{aligned}
\end{equation}
In this case, we say that $Q$ admits a $t$-admissible pair for $(r(\cdot),w)$.
\end{definition}

\begin{lemma}\label{lem: t adm modular inequality}
Let $1<r(\cdot)<\infty$ be such that $r_+<\infty$ and $[w]_{A_{r(\cdot)}}<\infty$.
Let $Q$ admit a $t$-admissible
pair for $(r(\cdot),w)$. Then for all $F\geq 0$ on $Q$
there holds that 
\[
\int_Q(F(x)w(x))^{r(x)}\,\mathrm d x\le t\Longrightarrow \int_Q\bigl(w(x)\avg{F}{Q}\bigr)^{r(x)}\,\mathrm d x
\lesssim_{r_+,[w]_{A_{r(\cdot)}},L}t.
\]
\end{lemma}

\begin{proof}
Let $(a,b)$ be a $t$-admissible pair for $(r(\cdot),w)$.
Young's inequality applied to
$F(x)w(x)$ and $b/w(x)$, together with (ii) and (iii) of
\eqref{t-admissible}, gives
\[
b|Q|\avg{F}{Q}
\le\int_Q(Fw)^{r(x)}\,\mathrm d x
+\int_Q(b/w)^{r'(x)}\,\mathrm d x
\le(1+L)t \leq 2[w]_{A_{r(\cdot)}}(1+L)ab|Q|.
\]
Thus,
$\avg{F}{Q}\le2[w]_{A_{r(\cdot)}}(1+L)a$. Thus, by (i),
\[
\int_Q\bigl(w(x)\avg{F}{Q}\bigr)^{r(x)}\,\mathrm d x
\le \bigl(2[w]_{A_{r(\cdot)}}(1+L)\bigr)^{r_+} Lt\lesssim_{r_+,[w]_{A_{r(\cdot)}},L}t.\qedhere
\]
\end{proof}

Admissible pairs can be interpolated as follows. 

\begin{lemma}\label{lem: t adm geo mean}Let $1<r(\cdot)<\infty$. Suppose that $(a_i,b_i)$ is $t_i$-admissible on $Q$ for $(r(\cdot),w)$,
$i=0,1$. Then, for all $\eta\in[0,1]$, the pair
$(a_0^\eta a_1^{1-\eta},b_0^\eta b_1^{1-\eta})$ is
$t_0^\eta t_1^{1-\eta}$-admissible on $Q$ for $(r(\cdot),w)$.
\end{lemma}

\begin{proof}
For a constant $0<\eta<1$, H\"older's inequality with exponents $1/\eta$ and
$1/(1-\eta)$ gives
\begin{align*}
\int_Q\bigl(w a_0^\eta a_1^{1-\eta}\bigr)^{r(x)}\,\mathrm d x
&=\int_Q
\bigl((w a_0)^{r(x)}\bigr)^\eta
\bigl((w a_1)^{r(x)}\bigr)^{1-\eta}\,\mathrm d x\\
&\le
\biggl(\int_Q(w a_0)^{r(x)}\,\mathrm d x\biggr)^\eta
\biggl(\int_Q(w a_1)^{r(x)}\,\mathrm d x\biggr)^{1-\eta}\le Lt_0^\eta t_1^{1-\eta}.
\end{align*}
This proves (i) of \eqref{t-admissible}. The same argument proves (ii),
and multiplying the corresponding powers of the two lower bounds
proves (iii). The cases $\eta=0$ and $\eta=1$ follow directly from
the hypotheses.
\end{proof}

The compatible
density gives an $\epsilon_Q$-admissible pair, for some $\epsilon_Q\leq 1$ (after rescaling $h$; see \eqref{eq:admPair0property} below), and $[w]_{A_{r(\cdot)}}<\infty$ gives a $1$-admissible pair.   
These pairs will be interpolated by Lemma \ref{lem: t adm geo mean} to prove Lemma \ref{lem:KEYlemma}. 

\begin{lemma}\label{lem:pairs2}
Let $1<p(\cdot)<\infty$ be such that $[w]_{A_{p(\cdot)}}<\infty$ and $h$ is a compatible density for $(p(\cdot),w).$ Let $Q$ be a fixed cube. 
\begin{itemize}
    \item The pair $(a_0,b_0) := \left(\left\langle\frac{wh^{1/p'(\cdot)}}{L}\right\rangle_Q,
\left\langle\frac{w^{-1}h^{1/p(\cdot)}}{L}\right\rangle_Q\right)$ is $\epsilon_Q$-admissible on $Q$ for $(\ppr,w^{-1})$, where 
\begin{align}\label{eq:defn:epsilonQ}
\epsilon_Q:=|Q|
\left\langle\frac{wh^{1/p'(\cdot)}}{L}\right\rangle_Q
\left\langle\frac{w^{-1}h^{1/p(\cdot)}}{L}\right\rangle_Q.
\end{align}
Moreover, if $E\subset Q$ is such that $|E|\geq \frac{1}{2}|Q|$ and  $\int h \leq \frac{1}{2}$,  then 
\begin{align}\label{eq:admPair0property}
    \epsilon_Q \leq \min\biggl\{\frac12,\,2\int_Eh\,\mathrm d x\biggr\} \leq 1.
\end{align}
\item The pair $(a_1,b_1) := \bigl(\|\mathbf 1_Q\|_{p'(\cdot),w^{-1}}^{-1},
\|\mathbf 1_Q\|_{p(\cdot),w}^{-1}\bigr)$ is $1$-admissible on $Q$ for $(\ppr,w^{-1})$.
\end{itemize}
\end{lemma}
\begin{proof}
We first verify \eqref{eq:admPair0property}.
By $h$ being a compatible density for $(p(\cdot),w)$, we have
$a_0 \leq w(x)h(x)^{1/p'(x)}$ and $b_0 \leq w(x)^{-1}h(x)^{1/p(x)}$ and thus 
\begin{align*} 
   \frac{\epsilon_Q}{|Q|} &:= a_0b_0 \leq \bigl(w(x)h(x)^{1/p'(x)}\bigr)
\bigl(w(x)^{-1}h(x)^{1/p(x)}\bigr) = h(x). 
\end{align*}
Integrating this over $Q$ and $E$ we have 
\[
\epsilon_Q\le\int_Qh\,\mathrm d x\le\frac12
\qquad\text{and}\qquad
\epsilon_Q\le2\int_Eh\,\mathrm d x
\] 
and thus \eqref{eq:admPair0property} follows.

To verify admissibility of the pair $(a_0,b_0)$, we need to check the three conditions in Definition \ref{defn:t-admissibility}. The bound on $a_0$ gives
\begin{multline*}
\int_Q\bigl(a_0/w(x)\bigr)^{p'(x)}\,\mathrm d x
=a_0\int_Qa_0^{p'(x)-1}w(x)^{-p'(x)}\,\mathrm d x\\ \le a_0\int_Q
\bigl(w(x)h(x)^{1/p'(x)}\bigr)^{p'(x)-1}
w(x)^{-p'(x)}\,\mathrm d x 
=a_0\int_Qw(x)^{-1}h(x)^{1/p(x)}\,\mathrm d x
=L\epsilon_Q,
\end{multline*}
and verifies (i). Interchanging $p(\cdot)$ and $p'(\cdot)$,
and $w$ and $w^{-1}$ verifies (ii). To verify (iii), by H\"older's inequality (Lemma \ref{lem:Holder}) we have 
\[
1 = \langle w w^{-1}\rangle_Q \leq 2\frac{\|\mathbf 1_Qw\|_{p(\cdot)}\|\mathbf 1_Qw^{-1}\|_{p'(\cdot)}}{|Q|} \leq 2[w]_{A_{p(\cdot)}}
\]
and thus 
\[
a_0b_0 =: \frac{\epsilon_Q}{|Q|} \geq \frac{\epsilon_Q}{2[w]_{A_{p(\cdot)}}|Q|}.
\]

Then we verify the $1$-admissibility of the pair $(a_1,b_1)$. By \eqref{eq:prelim-indicator-normalization} applied
to $(p(\cdot),w)$ and $(p'(\cdot),w^{-1})$ we have 
\[
\int_Q\biggl(\frac{w(x)}{\|\mathbf 1_Q\|_{p(\cdot),w}}\biggr)^{p(x)}
\,\mathrm d x
=\int_Q\biggl(\frac{w(x)^{-1}}{\|\mathbf 1_Q\|_{p'(\cdot),w^{-1}}}\biggr)^{p'(x)}
\,\mathrm d x
=1\le L.
\]
These are (i) and (ii). That $[w]_{A_{p(\cdot)}}<\infty$ gives (iii).
\end{proof}

The following lemma allows rescaling the compatible density so that \eqref{eq:admPair0property} can be used. 
\begin{lemma}\label{lem:rescaleH}
If $h$ is a compatible density for $(\pp,w)$, where $1<p(\cdot)<\infty$, then so is $ch$, for all $c>0.$
\end{lemma}
\begin{proof} If $c\geq 1$, then $c^{1/p(\cdot)}\leq c$ and if $c\leq 1$, then $c^{1/p(\cdot)}\leq 1;$ thus, $c^{1/p(\cdot)} \leq \max(1,c).$ Thus,  
\begin{multline}
    \langle w^{-1}(ch)^{1/p(\cdot)} \rangle_Q \leq \max(1,c) \langle w^{-1}h^{1/p(\cdot)}\rangle_Q  \\
\leq \max(1,c)\left[ w^{-1}h^{1/p(\cdot)}\right]_{A_1}\essinf_Qw^{-1}h^{1/p(\cdot)} \\ \lesssim \max(1,c)\max(1,c^{-1})\left[ w^{-1}h^{1/p(\cdot)}\right]_{A_1}\essinf_Qw^{-1}(ch)^{1/p(\cdot)}.
\end{multline}   
The same bound holds with $p'(\cdot)$ in place of $p(\cdot)$ and $w$ in place of $w^{-1}$.
\end{proof}

\begin{proof}[Proof of Lemma~\ref{lem:KEYlemma}]
Let $p(\cdot)$ admit a compatible density $h$ and satisfy $[w]_{A_{p(\cdot)}}<\infty$.  By Lemma \ref{lem:rescaleH}, substituting $\frac{1}{2}\|h\|_{L^1}^{-1}\cdot h$ for $h$, we may assume that $\int h\leq \frac{1}{2}.$
Let $(a_0,b_0)$ be the $\epsilon_Q$-admissible and $(a_1,b_1)$ the
$1$-admissible pairs of Lemma~\ref{lem:pairs2}.
 By $|E|\geq \frac{1}{2}|Q|$ (standing assumption) and  $\int h\leq \frac{1}{2}$, it follows that \eqref{eq:admPair0property} is in force, i.e. that $\epsilon_Q \leq  1$. Together with  $\int_E(f(x)/w(x))^{p'(x)}\,\mathrm d x\le1$ (standing assumption) this yields
\[
\epsilon_Q\leq t:=\max\biggl\{
\int_E\Bigl(\frac{f(x)}{w(x)}\Bigr)^{p'(x)}\,\mathrm d x,
\epsilon_Q\biggr\} \leq 1.
\]
Now let $\eta\in[0,1]$ be such that
$t=\epsilon_Q^\eta$. By Lemma~\ref{lem: t adm geo mean} there exists a
$t$-admissible pair on $Q$ for $(p'(\cdot),w^{-1})$.
Hence, by Lemma~\ref{lem: t adm modular inequality}, applied with $F=f1_E$, we obtain
\begin{align*}
\int_Q\biggl(\frac{\avg{f1_E}{Q}}{w(x)}\biggr)^{p'(x)}\,\mathrm d x
\lesssim_{(p')_+,[w]_{A_{p(\cdot)}},L}t
\le  \int_E\Bigl(\frac{f(x)}{w(x)}\Bigr)^{p'(x)}\,\mathrm d x + \epsilon_Q \leq \int_E\Bigl(\frac{f(x)}{w(x)}\Bigr)^{p'(x)}\,\mathrm d x
+2\int_Eh(x)\,\mathrm d x,
\end{align*}
where we used \eqref{eq:admPair0property} in the last step.
\end{proof}

\subsection{Necessity}\label{sec:proof2}
Assuming that $M$ is bounded on $L^{p(\cdot)}_w$, we prove that $[w]_{A_{p(\cdot)}}<\infty$ and construct a compatible density $h$ for $(p(\cdot),w)$. 
The Muckenhoupt type condition follows from a duality argument used in the unweighted setting by Kopaliani \cite{Kop07}; see also \cite[Corollary~4.50]{CUF13}. In fact, this necessity holds far more generally: Berezhnoi \cite{Ber99} showed that if $M$ is bounded on an ideal Banach space $X$, then $\sup_Q |Q|^{-1}\|\mathbf 1_Q\|_X\|\mathbf 1_Q\|_{X'}<\infty$. We include the short proof of this fact for $X = L^{\pp}_w$ for completeness. To construct the compatible density, we use a nonlinear iteration based on the Rubio de Francia algorithm (see Rubio de Francia \cite[Theorem~1 and Section~3]{RdF82} and Coifman, Jones and Rubio de Francia \cite[proof of Theorem~II]{CJR83}). A related construction appears in Rutsky \cite[proof of Theorem~8]{Rut14}, who seeks a majorant in $A_q$ whose $\delta$-th power belongs to $A_1$, for fixed $q>1$ and $\delta>0$.

\begin{proof}[Proof of Theorem \ref{thm:main}; necessity]
For a fixed cube $Q$, set $E_m:=Q\cap\{m^{-1}\le w\le m\}$, $m\ge1$. Since $w$ and $w^{-1}$ are bounded on $E_m$, both indicator norms below are finite. 
Applying Lemma~\ref{lem:dualityFormula} to compute $\|1_{E_m}\|_{p'(\cdot),w^{-1}}$ we obtain
\begin{align*}
\frac{\|1_{E_m}\|_{\pp,w}
      \|1_{E_m}\|_{\ppr,w^{-1}}}{|Q|}
&\simeq
\sup_{\|f\|_{\pp,w}\le1}
\bigl\|\avg{|f|1_{E_m}}{Q}1_{E_m}\bigr\|_{\pp,w}\\
&\le
\sup_{\|f\|_{\pp,w}\le1}
\|M(f1_{E_m})\|_{\pp,w}
\le \|M\|_{L^{\pp}_w\to L^{\pp}_w}.
\end{align*}
As $E_m\uparrow Q$, Proposition~\ref{prop:Lp=Fatou} allows us to let $m\to\infty$.
Taking supremum over cubes yields $[w]_{A_{p(\cdot)}}<\infty$.

Since $w \in A_{\pp}$ and $1 < p_- \leq p_+ < \infty$, Remark~\ref{remark:LocalIntegrability} implies $w^{\pp}, w^{-\ppr} \in L^1_{\mathrm{loc}}$. So, by Theorem~\ref{thm:prelim-maximal-duality}, $M$ is bounded on $L^{p'(\cdot)}_{w^{-1}}$ and hence 
\[
N:=\max\bigl\{
\|M\|_{L^{p(\cdot)}_w\to L^{p(\cdot)}_w},
\|M\|_{L^{p'(\cdot)}_{w^{-1}}\to 
L^{\ppr}_{w^{-1}}}
\bigr\} < \infty.
\]
Choose a positive function $\varphi\in L^1(\R^n)$ with $\int_{\R^n}\varphi\,\mathrm d x=1/3$. Starting from $H_0:=\varphi$, define recursively
\begin{align}\label{eq:PHI}
    H_{k+1}:=\max\Biggl\{\varphi,
\Biggl(\frac{wM\bigl(w^{-1}H_k^{1/\pp}\bigr)}{3N}\Biggr)^{\pp},\Biggl(\frac{w^{-1}M\bigl(wH_k^{1/\ppr}\bigr)}{3N}\Biggr)^{\ppr}\Biggr\} =: \Phi(H_{k}).
\end{align}

The normalization by $3N$ ensures that $\int_{\R^n}H_k\,\mathrm d x\leq1$ for all $k\geq0$, which we will now prove  by induction.
The case $k=0$ follows from the choice of $\varphi$. Suppose $\int_{\R^n}H_k\,\mathrm d x\leq1$. For either choice $s(\cdot)=p(\cdot)$ and $u=w$, or $s(\cdot)=p'(\cdot)$ and $u=w^{-1}$, we have
\begin{align}\label{eq:XX}
    \rho_{s,u}\bigl(u^{-1}H_k^{1/s(\cdot)}\bigr)
=\int_{\R^n}H_k\,\mathrm d x\leq1.
\end{align}
By \eqref{eq:prelim-unit-ball} of Lemma~\ref{lem:prelim-norm-modular}, this gives $\|u^{-1}H_k^{1/s(\cdot)}\|_{s(\cdot),u}\leq1$. Thus $\|M(u^{-1}H_k^{1/s(\cdot)})/N\|_{s(\cdot),u} \leq 1$, by the definition of $N$, and hence, applying \ref{eq:prelim-unit-ball} once more, gives
$\rho_{s(\cdot),u}\bigl(M(u^{-1}H_k^{1/s(\cdot)})/N\bigr)\leq1$.
Since $s(x)\ge1$ almost everywhere, it follows that
\begin{equation}\label{eq:stepp4}
\rho_{s,u}\Biggl(\frac{M\bigl(u^{-1}H_k^{1/s(\cdot)}\bigr)}{3N}\Biggr)
\leq\frac13\rho_{s,u}\Biggl(\frac{M\bigl(u^{-1}H_k^{1/s(\cdot)}\bigr)}{N}\Biggr)
\leq\frac13.
\end{equation}
Bounding the maximum by the sum and using \eqref{eq:stepp4} gives
\begin{multline*}\label{eq:stepp1}
\int_{\R^n}H_{k+1}\,\mathrm d x
\leq\int_{\R^n}\varphi\,\mathrm d x
+\rho_{p(\cdot),w}\Biggl(\frac{M\bigl(w^{-1}H_k^{1/\pp}\bigr)}{3N}\Biggr)+\rho_{p'(\cdot),w^{-1}}\Biggl(\frac{M\bigl(wH_k^{1/\ppr}\bigr)}{3N}\Biggr)
\leq\frac13+\frac13+\frac13=1,
\end{multline*}
which completes the induction.

Then we verify that $H_k\le H_{k+1}$. Clearly $H_0=\varphi\leq H_1$, for the base case. 
For the inductive step we use that if $0\leq f\leq g,$ then $Mf\leq Mg$, and this gives $\Phi(f)\leq \Phi(g),$ where $\Phi$ is as on the line \eqref{eq:PHI}. Thus, using the inductive assumption $H_{k-1}\le H_{k}$, we have 
$H_{k} := \Phi(H_{k-1}) \leq \Phi(H_k) =: H_{k+1}.$

Since $H_k$ increases, the pointwise limit
\begin{equation}\label{eq:Seed->CompDens}
H:=\lim_{k\to\infty}H_k
\end{equation}
exists in $[0,\infty]$. By monotone convergence and
\eqref{eq:XX},
\begin{equation*}\label{eq:Seed->CompDens2}
\int_{\R^n}H\,\mathrm d x
=\lim_{k\to\infty}\int_{\R^n}H_k\,\mathrm d x
\le1.
\end{equation*}
Thus $H<\infty$ almost everywhere and $H\in L^1$.

Both $w^{-1}H^{1/\pp}$ and $wH^{1/\ppr}$ are positive,
since $H\ge H_0=\varphi>0$. To prove that they belong to
$A_1$, we first check local integrability.
For either choice $s(\cdot)=p(\cdot)$ and $u=w$, or $s(\cdot)=p'(\cdot)$ and $u=w^{-1}$, we have 
$\rho_{s(\cdot),u}(u^{-1}H^{1/s(\cdot)}) = \int_{\R^n}H \leq 1.$ Thus, by \eqref{eq:prelim-unit-ball},  $\|u^{-1}H^{1/s(\cdot)}\|_{s(\cdot), u} \leq 1$ and hence by $L^{s(\cdot)}_u\subset L^1_{\mathrm{loc}}$ (see Remark \ref{remark:LocalIntegrability}) that  $u^{-1}H^{1/s(\cdot)} \in L^1_{\mathrm{loc}}.$
Thus monotone convergence and \eqref{eq:PHI} give, for a fixed cube $Q$ and almost every $x\in Q$, that 
\begin{multline}
    \bigl\langle u^{-1}H^{1/s(\cdot)}\bigr\rangle_Q
=\lim_{k\to\infty}
  \bigl\langle u^{-1}H_k^{1/s(\cdot)}\bigr\rangle_Q \leq \lim_{k\to\infty} 
   M\left(u^{-1}H_k^{1/s(\cdot)}\right)(x) \\ 
\leq \lim_{k\to\infty} 3N u(x)^{-1}H_{k+1}^{1/s(x)}=  3N u(x)^{-1}H(x)^{1/s(x)}.
\end{multline}
Thus $w^{-1}H^{1/\pp}, wH^{1/\ppr}\in A_1$, by \eqref{eq:defnA1}.
\end{proof}

\section{Applications and examples}\label{sec:Applications}
In this section we exemplify some applications of Theorem \ref{thm:main} by reproving and in some cases extending, known theorems using the framework of compatible densities. As we will show, this produces rather short and clean proofs. It is not exhaustive, but rather a proof of concept. 

We begin with giving a new proof that $\pp\in LH_\infty$ and $1\in A_{p({\cdot})}$ imply $M : L^{\pp} \to L^{\pp}$. This is known; it follows from the results of \cite[Theorem 4.52]{CUF13}, where $LH_\infty$ is replaced by the weaker condition $N_\infty$ of Nekvinda \cite{Nek04}. (We recall that an alternative notation for $1\in A_{p({\cdot})}$ is $p(\cdot)\in K_0$.)
\begin{definition}
    Given an exponent function $\pp$, we say $\pp \in LH_\infty$ if there exist constants $C > 0$ and $p_\infty \in [p_-, p_+]$ such that for a.e. $x \in \R^n$,
    \begin{equation} \label{defn:log Holder}
        \left| \frac{1}{p(x)} - \frac{1}{p_\infty} \right| \leq \frac{C}{\log{(e + |x|)}}.
    \end{equation}
\end{definition}

The idea is that $LH_\infty$ alone allows you to build a compatible density via power weights.  
\begin{theorem} \label{thm:log holder}
    Let $1 < p_- \leq p_+ < \infty$. If $p(\cdot)\in LH_\infty$ and $1\in A_{p(\cdot)}$, then $M : L^{\pp} \to L^{\pp}$ boundedly.
\end{theorem}
\begin{proof}
    Let $\pp \in LH_\infty$ and $C,p_\infty$ be as in \eqref{defn:log Holder}. Setting $p'_\infty := (p_\infty)'$, we see $\ppr \in LH_\infty$ as well since
    \begin{equation} \label{eqn: LHinfty}
        \left| \frac{1}{p'(x)} - \frac{1}{p_\infty'} \right| = \left| \frac{1}{p(x)} - \frac{1}{p_\infty} \right| \leq \frac{C}{\log{(e + |x|)}}.
    \end{equation}   
    Pick $N$ satisfying $n < N < n\min\{p_\infty,p_\infty'\}$, and define $h(x) := \min\{1, |x|^{-N}\}$. We claim that $h$ is a compatible density for $\pp$. 

    It is clear that $h \in L^1$ since $N > n$. For $\qq \in \{\pp,\ppr\}$, notice $h^{1/\qq} = \min\{1,|x|^{-N/\qq}\}$. Thus for $|x| > 1$, the assumption $\qq \in LH_\infty$ implies
   \[ 
    \frac{h^{1/q(x)}}{|x|^{-N/q_\infty}} = |x|^{-N\left(\tfrac1{q(x)} - \tfrac1{q_\infty}\right)} = e^{-t(x)}, \qquad t(x) := N\left(\frac{1}{q(x)} - \frac{1}{q_\infty}\right) \log|x|. \]
    Then by \eqref{eqn: LHinfty},
    \[ 
    |t(x)| \leq CN\frac{\log{|x|}}{\log{(e+|x|)}} \leq CN, 
    \]
    and thus $\frac{h^{1/q(x)}}{|x|^{-N/q_\infty}} \gtrsim 1$. Bounding in the same manner, we have $\frac{|x|^{-N/q_\infty}}{h^{1/q(x)}} \gtrsim 1$. It follows that 
    $$
    h^{1/\qq} \simeq \min\{1,|x|^{-N/q_\infty}\}.
    $$ 
    It is standard that if $0 \leq \alpha < n$, then $|x|^{-\alpha} \in A_1$. Therefore, since $0 < N/q_\infty < n$, $|x|^{-N/q_\infty} \in A_1$. Recalling that $A_1$ is stable under taking minimums and functions comparable to $A_1$ weights are themselves $A_1$ weights, we conclude $h^{1/\qq} \in A_1$. 
\end{proof}
We can also recover a result of Lerner \cite[Corollary 6.5]{LerCrit26}. Here we offer a shorter proof built on lighter machinery and remove the $1/n$ factor in the condition on $\alpha$. 
\begin{theorem} \label{thm:radial}
    Suppose that $p(\cdot)$ is a radial exponent such that $1 < p_- \leq p_+ < \infty$ and $1 \in A_{p(\cdot)}.$ If there exist $N > 1$ and $\alpha$ satisfying 
    \begin{align}\label{eqn:lerner reproof1}
        0 < \alpha < p_- \min\{1,(p_- - 1)\} 
    \end{align}
    such that for all $|x| \geq |y| \geq N$, 
    \begin{equation} \label{eqn:lerner reproof}
        |p(x) - p(y)| \leq \alpha \frac{1}{\log{|y|}} \log{\frac{|x|}{|y|}},
    \end{equation}  
    then $M : L^{\pp} \to L^{\pp}$. 
\end{theorem}
We will use the following elementary lemma about radial $A_1$ weights. It is essentially well known; we include its short proof for completeness. 
\begin{lemma} \label{lem:radial A1}
    Let $w(x) = W(|x|)$ with $W > 0$ bounded above and below on $[0,R_0]$. Suppose that on $[R_0, \infty)$ the function $W$ is nonincreasing and $r^\gamma W(r)$ is nondecreasing for some $\gamma < n$. Then $w \in A_1$.   
\end{lemma}
\begin{proof}
    Since $W$ is bounded from above and below on $[0,R_0]$ and is nonincreasing on $[R_0,\infty)$, it follows that $W$ is ``almost nonincreasing on $[0,\infty)$", i.e., $W(r) \lesssim W(s)$ for $s \leq r$. Moreover, $r^\gamma W(r)$ being nondecreasing on $[R_0,\infty)$ implies
    \[ \left(\frac{r}{s}\right)^{-\gamma} \leq \frac{W(r)}{W(s)} \leq 1, \qquad R_0 \leq s \leq r.\]
    This, along with the assumption that $W$ is bounded from below and above on $[0,R_0]$, gives that $W$ is doubling: $W(r) \simeq W(s)$, for $s \leq r \leq 2s$. 

    Fix $x \in \R^n$ and a cube $Q$ containing $x$ with side length $\ell$. Suppose that $\sqrt{n} \ell \leq |x|/2$. Then $Q \subset \{y : |x|/2 \leq |y| \leq 2|x|\}$. By our estimates above, this implies $w(y) \simeq w(x)$ for all $y \in Q$, and so $\langle w \rangle_Q \simeq w(x)$. Now suppose $\sqrt{n} \ell > |x|/2$. Then $Q \subset B := B(0,\rho)$ where $\rho = 3\sqrt{n}\ell$. If $\rho \leq R_0$ then $w$ is bounded above and below on $Q$, so trivially $\langle w \rangle_Q \lesssim w(x)$. So, assume $\rho > R_0$. Since $|Q| \simeq |B(0,\rho)|$ and $w > 0$, we have
    %
    \[ \langle w \rangle_Q \lesssim  \langle w \rangle_B = \frac{1}{|B|} \int_B w(x) \, \mathrm d x \lesssim \rho^{-n} \int_0^\rho W(r) r^{n-1} \, dr = \rho^{-n} \left(\int_0^{R_0} + \int_{R_0}^{\rho}\right) W(r) r^{n-1} \, dr. \]
  Since $W$ is bounded on $[0,R_0]$, $\rho \geq R_0, \gamma < n$, and $r^{\gamma}W(r)$ is nondecreasing on $[R_0,\infty)$,
    \[ \int_0^{R_0} W(r)r^{n-1} \, dr \lesssim 1 \sim W(R_0) R_0^n \leq W(R_0) R_0^\gamma \rho^{n-\gamma} \leq W(\rho) \rho^\gamma \rho^{n-\gamma} = W(\rho) \rho^n. \]
      Since $r^\gamma W(r)$ is nondecreasing on $[R_0,\infty)$ and $\gamma < n$,
    \[ \int_{R_0}^\rho W(r) r^{n-1} \, dr \leq W(\rho) \rho^{\gamma} \int_{R_0}^\rho r^{n-\gamma-1} \, dr \lesssim W(\rho) \rho^n. \]
    Hence $
    \langle w \rangle_Q \lesssim W(\rho)
    $.
    Because $W$ is almost nonincreasing, $W(r) \lesssim W(s)$, for $s \leq r$, it follows that $\langle w \rangle_Q \lesssim w(x)$, for all $x\in Q$. This confirms that $w \in A_1$. 
\end{proof}
\begin{proof}[Proof of Theorem \ref{thm:radial}]
     Since we assume $1\in A_{p(\cdot)}$, by Theorem \ref{thm:main} it suffices to show $\eqref{eqn:lerner reproof}$ provides a compatible density for $\pp$. We will show that $h(x) := \max\{|x|,N\}^{-\beta}$ is a compatible density for all $\beta \in (n,n/\kappa)$, where
     \[ \kappa := \max\Big\{\frac{p_-+\alpha}{p_-^2}, \, 1-\frac{p_--\alpha}{p_-p_+}\Big\}<1 .\]
     Indeed, by \eqref{eqn:lerner reproof1} the first item inside the maximum is strictly bounded from above by $1.$
     It is clear that $h \in L^1$ since $\beta > n$. 
     Let $\qq \in \{\pp, \ppr\}$ and define 
     $$
     W = h^{1/\qq} = \max\{|x|,N\}^{-\beta/\qq}.
     $$
     We show that $W\in A_1$ by verifying that $W$ satisfies the hypotheses of Lemma \ref{lem:radial A1}.  
     
     On $[0,N]$, $W$ is bounded above and below. On $[N,\infty)$ we have $W=|x|^{-\beta/\qq}.$ To establish the nonincreasing and nondecreasing properties, we first establish some pointwise estimates on $W$. For $N \leq s \leq r$ we have 
     \[ -\frac{1}{\beta} \log{\frac{W(r)}{W(s)}} = \frac{1}{q(r)} \log{r} - \frac{1}{q(s)}\log{s} = \frac{1}{q(r)} \log{\frac{r}{s}} + \left(\frac{1}{q(r)} - \frac{1}{q(s)}\right)\log{s}. \]
     An exponent computation, and using the bound \eqref{eqn:lerner reproof}, shows that for both $\qq = \pp,\ppr$,
     \[ 
     \left|\frac{1}{q(r)} - \frac{1}{q(s)}\right| = \frac{|p(s) - p(r)|}{p(r)p(s)} \leq \frac{\alpha}{p(r)p(s)}\frac{1}{\log(s)} \log{\frac{r}{s}}. 
     \]
     This gives
     \[ \left(\frac{1}{q(r)} - \frac{\alpha}{p(r)p(s)} \right) \log{\frac{r}{s}} \leq -\frac{1}{\beta} \log{\frac{W(r)}{W(s)}} \leq \left(\frac{1}{q(r)} + \frac{\alpha}{p(r)p(s)} \right) \log{\frac{r}{s}}. \]
     Multiplying by $-\beta$ and exponentiating gives
     \begin{equation} \label{eqn:radial est}
         \left(\frac{r}{s}\right)^{-\beta\left(\tfrac{1}{q(r)} + \tfrac{\alpha}{p(r)p(s)} \right)} \leq \frac{W(r)}{W(s)} \leq \left(\frac{r}{s}\right)^{-\beta\left(\tfrac{1}{q(r)} - \tfrac{\alpha}{p(r)p(s)} \right)}.
     \end{equation}
     
     To show that $W$ is nonincreasing on $[N,\infty)$, it suffices to show that the exponent on the right hand side of \eqref{eqn:radial est} is negative. We do this individually for $\pp$ and $\ppr$. By our assumptions on $\alpha$,
     \[ \frac{\alpha}{p(r)p(s)} \leq \frac{p_-}{p(r)p(s)} \leq \frac{1}{p(r)}, \]
     and
     \[ \frac{\alpha}{p(r)p(s)} \leq \frac{p_-(p_- - 1)}{p(r)p(s)} \leq \frac{p(s)(p(r) - 1)}{p(r)p(s)} = \frac{1}{p'(r)}. \]

     To conclude the proof, we show that $r^\gamma W(r)$ is nondecreasing for $\gamma = \beta \kappa < n$.  By the lower bound in \eqref{eqn:radial est}, it suffices to show
     \[ \gamma \geq \beta\left(\frac{1}{q(r)} + \frac{\alpha}{p(r)p(s)}\right), \]
     which by our definition of $\gamma$, is equivalent to showing
     \[ \kappa \geq \frac{1}{q(r)} + \frac{\alpha}{p(r)p(s)}. \]
      We verify this individually for $\pp$ and $\ppr$.
     By our definition of $\kappa$, 
     \[ \frac{1}{p(r)} + \frac{\alpha}{p(r)p(s)} \leq \frac{1}{p_-} + \frac{\alpha}{p_-^2} = \frac{p_- + \alpha}{p_-^2} \leq \kappa, \]
     and since $\alpha \leq p_-$, 
     \[ \frac{1}{p'(r)} + \frac{\alpha}{p(r)p(s)} = 1 - \frac{1}{p(r)} \left(1 - \frac{\alpha}{p(s)}\right) \leq 1 - \frac{1}{p_+} \left(1 - \frac{\alpha}{p_-}\right) = 1 - \frac{p_- - \alpha}{p_-p_+} \leq \kappa. \]
\end{proof}
\section{Failure of the Diening condition for \texorpdfstring{$L^{p(\cdot)}_w$}{Lp(.) (w)}}
\label{sec:failure-diening}

Diening \cite[Theorem~8.1]{Die05} proved that, for $1<p_-\leq p_+<\infty$, boundedness of $M$ on $L^{\pp}$ is equivalent to uniform boundedness of the operators
\[
T_{\mathcal F}f:=\sum_{Q\in\mathcal F}\avg{|f|}{Q}1_Q,
\]
where $\mathcal F$ ranges over finite families of pairwise disjoint cubes. Uniform boundedness of the $T_{\mathcal F}$ on a Banach function space $X$ is referred to as $X$ satisfying the strong Muckenhoupt condition (see \cite[Definition~3.5]{Nie26}).
As we show in our next proposition, the strong Muckenhoupt condition does not characterize boundedness of $M$ on the weighted variable Lebesgue spaces.
Although the condition remains necessary on $L^{\pp}_w$ (since $T_{\mathcal F}f\leq Mf$), the following example shows that it is no longer sufficient. 
This failure further motivates our characterization, Theorem \ref{thm:main}, which remains valid in the weighted setting. 
It also gives a negative answer to the additional question immediately following Conjecture~D in \cite{Nie26}.

\begin{proposition}\label{prop:dc-counterexample}
There exist an exponent function $\pp$ and a weight $w$ on $\R$, with $1<p_-\leq p_+<\infty$ and $w^{\pp},w^{-\ppr}\in L^1_{\mathrm{loc}}(\R)$, such that
\[
\sup_{\mathcal F}
\|T_{\mathcal F}\|_{L^{p(\cdot)}_w\to L^{p(\cdot)}_w}<\infty,
\]
but $M$ is unbounded on $L^{p(\cdot)}_w$.
\end{proposition}

\begin{proof}
Define
\[
p(x)=\begin{cases}2,&x\leq0,\\4,&x>0,\end{cases}
\qquad
w(x)=\begin{cases}1,&x\leq0,\\x^{1/4},&x>0.\end{cases}
\]
The local integrability assumptions hold for $w$ since $w^{\pp}=x$ and $w^{-\ppr}=x^{-1/3}$ on $(0,\infty)$. 
Then
\begin{equation} \label{eqn:modular identity}
\rho_{\pp,w}(f)=\int_{-\infty}^0|f(x)|^2\,\mathrm dx
+\int_0^\infty|f(x)|^4x\,\mathrm dx.
\end{equation}

Fix a finite disjoint family $\mathcal F$ and $f\geq0$ with $\rho_{\pp,w}(f)\leq1$. 
By homogeneity and Lemma~\ref{lem:prelim-norm-modular}, it suffices to prove $\rho_{\pp,w}(T_{\mathcal F}f)\lesssim1$, uniformly in $f$ and $\mathcal F$.
For $J=(a,b)\subset(0,\infty)$,
\[
\avg{x}{J}\leq b,
\qquad
\avg{x^{-1/3}}{J}
=\frac{3(b^{2/3}-a^{2/3})}{2(b-a)}
\leq\frac32b^{-1/3}.
\]
Hence, H\"older's inequality gives
\[ \rho_{\pp,w}(\langle f \rangle_J 1_J) = \langle f \rangle_J^4 \int_J x \, dx \leq \frac{1}{|J|^4} \int_J f(x)^4 x \, dx \left(\int_J x^{-1/3} \, dx \right)^3 \left(\int_J x \, dx\right) \lesssim \rho_{\pp,w}(f 1_J). \]
The same estimate for $J\subset(-\infty,0)$ follows from
Jensen's inequality;
\[ 
\rho_{\pp,w}(\langle f \rangle_J 1_J) = \int_J \langle f \rangle_J^2 \, dx \leq \int_J f(x)^2 \, dx = \rho_{\pp,w}(f 1_J). 
\]
For $I=(-a,b)$ with $a,b>0$, we split 
\[ \langle f \rangle_I = \frac{1}{|I|} \left(\int_{-a}^0 f(x) \, dx + \int_0^b f(x) \, dx\right). \]
The first integral can be estimated by H\"older's inequality;
\[ \int_{-a}^0 f(x) \, dx \leq a^{1/2} \left(\int_{-a}^0 f(x)^2 \, dx\right)^{1/2} \leq |I|^{1/2} \rho_{\pp,w}(f 1_{(-a,0)})^{1/2} \leq |I|^{1/2}. \]
Similarly, for the second integral,
\[ \int_0^b f(x) \, dx \leq \left(\int_0^b f(x)^4 x \, dx\right)^{1/4} \left(\int_0^b x^{-1/3} \, dx \right)^{3/4} \lesssim |I|^{1/2}. \]
Together, this gives $\langle f \rangle_I \lesssim |I|^{-1/2}$. So, writing $I = (-a,0] \cup (0,b)$ and applying \eqref{eqn:modular identity}, 
\[ 
\rho_{\pp,w}(\langle f \rangle_I1_I) = a \langle f \rangle_I^2 + \frac{b^2}{2} \langle f \rangle_I^4 \lesssim \frac{a}{|I|} + \frac{b^2}{|I|^2} \leq 1. \]
Since a disjoint family $\mathcal F$ contains at most one interval
meeting both half-lines, these estimates give
\[
\rho_{p(\cdot),w}(T_{\mathcal F}f)
=\sum_{J\in\mathcal F}\rho_{p(\cdot),w}(\avg{f}{J}1_J)
\lesssim\rho_{p(\cdot),w}(f)+1\lesssim1.
\]

To prove $M$ is not bounded on $L^{p(\cdot)}_w$, set $f(x)=(x\log(e/x))^{-1/2} 1_{(0,1)}(x)$. Then
\[
\rho_{p(\cdot),w}(f)=\int_0^1\frac{ dx}{x\log^2(e/x)}=1.
\]
Since $f$ is decreasing on $(0,1)$, for $0<t<1$ we have
\[
Mf(-t)\geq\frac1{2t}\int_0^t f(x)\, dx
\geq\frac12f(t).
\]
Consequently,
\[
\int_{-1}^0 Mf(t)^2\, dt
\geq\frac14\int_0^1\frac{dx}{x\log(e/x)}=\infty,
\]
so $Mf\notin L^{p(\cdot)}_w$.
\end{proof}
\section*{Acknowledgments}
We wish to thank Emiel Lorist for telling us about the counterexample in Proposition \ref{prop:dc-counterexample}, which was found using an LLM.
\section*{AI disclosure statement}
An LLM was used to explore characterizations for the boundedness of $M$ in unweighted variable Lebesgue spaces. 
More precisely, in a long initial conversation with an LLM we iteratively constructed a very rough proof of Theorem \ref{thm:main-unweighted}. 
Notably this proof utilized Diening's condition and did not use sparse domination.
We then worked to simplify the proof to the unweighted version of the current proof.
Afterwards we used an LLM to obtain the proof of Theorem \ref{thm:main} from the proof of Theorem \ref{thm:main-unweighted}, which generalized smoothly.

An LLM was used throughout the process to hunt for errors, typos and relevant references.

All mathematical arguments in their final form were written and verified by the authors, who take full responsibility for all claims and arguments in the paper.

\bibliographystyle{alpha}
\bibliography{references}
\end{document}